\documentclass[12pt, a4paper]{article}
\usepackage{amsfonts}

\usepackage{graphicx}
\usepackage{amsmath}
\usepackage[french]{babel}

\newtheorem{theorem}{Th\'{e}or\`{e}me}[section]

\newtheorem{conjecture}[theorem]{Conjecture}
\newtheorem{corollary}[theorem]{Corollaire}

\newtheorem{definition}[theorem]{D\'{e}finition}

\newtheorem{lemma}[theorem]{Lemme}

\newtheorem{proposition}[theorem]{Proposition}
\newtheorem{remark}[theorem]{Remarque}

\newenvironment{proof}[1][D\'{e}monstration]{\noindent\textbf{#1.} }{\hfill $\Box $}

\begin{document}

\title{Etude des jets de Demailly-Semple en dimension 3}
\author{Erwan Rousseau}
\date{}
\maketitle

\begin{abstract}
Dans cet article nous faisons l'\'{e}tude alg\'{e}brique des jets de
Demailly-Semple en dimension 3 en utilisant la th\'{e}orie des invariants
des groupes non r\'{e}ductifs. Cette \'{e}tude fournit la
caract\'{e}risation g\'{e}om\'{e}trique du fibr\'{e} des jets d'ordre 3 sur
une vari\'{e}t\'{e} de dimension 3 et permet d'effectuer, par Riemann-Roch,
un calcul de caract\'{e}ristique d'Euler.
\end{abstract}

\section{Introduction}

\subsection{Contexte g\'{e}om\'{e}trique}

Il est bien connu (cf. \cite{De95}, \cite{SY95}) que l'\'{e}tude de
l'hyperbolicit\'{e} des vari\'{e}t\'{e}s alg\'{e}briques complexes est
li\'{e}e \`{a} l'\'{e}tude des sections globales de certains fibr\'{e}s
vectoriels $E_{k,m}T_{X}^{\ast }$ d'op\'{e}rateurs diff\'{e}rentiels d'ordre
$k$ et de degr\'{e} $m$ agissant sur les germes de courbes holomorphes dans $%
X$, vari\'{e}t\'{e} complexe.

L'\'{e}tude des jets de Demailly-Semple est motiv\'{e}e par les
r\'{e}sultats qu'ils ont fournis sur l'hyperbolicit\'{e} des
vari\'{e}t\'{e}s complexes, sur des questions li\'{e}es \`{a} la conjecture
de Kobayashi qui stipule que le compl\'{e}mentaire d'une hypersurface
g\'{e}n\'{e}rique de degr\'{e} $d\geq 2n+1$ dans $\mathbb{P}_{\mathbb{C}%
}^{n} $ est hyperbolique.

Des r\'{e}sultats int\'{e}ressants ont \'{e}t\'{e} obtenus en dimension 2
pour le cas des compl\'{e}mentaires de courbes dans $\mathbb{P}_{\mathbb{C}%
}^{2}$ avec un nombre donn\'{e} $k$ de composantes irr\'{e}ductibles: citons
les r\'{e}sultats de Y.T. Siu et S.K. Yeung \cite{SY95}, ceux de J.P.
Demailly et J. El Goul \cite{DEG00} qui ont trait\'{e} le cas du
compl\'{e}mentaire d'une courbe g\'{e}n\'{e}rique lisse dans $\mathbb{P}_{%
\mathbb{C}}^{2},$ ainsi que ceux de G. Dethloff, G. Schumacher et P.M. Wong
\cite{DSW92} et \cite{DSW94} concernant le cas de 3 composantes, et plus
r\'{e}cemment le cas de 2 composantes a \'{e}t\'{e} trait\'{e} \'{e}galement
par les techniques de jets \cite{Rou03}.

L'hyperbolicit\'{e} en dimension 3 est un sujet tr\`{e}s peu
d\'{e}frich\'{e} \`{a} l'heure actuelle malgr\'{e} l'existence de quelques
classes d'exemples (Masuda-Noguchi, Siu, Shiffmann et Zaidenberg).
L'\'{e}tude des jets de Demailly-Semple en dimension 3 a \'{e}t\'{e} faite
dans la perspective d'attaquer le probl\`{e}me de l'hyperbolicit\'{e} des
hypersurfaces projectives g\'{e}n\'{e}riques de grand degr\'{e} de dimension
3 pour lequel il n'y a pas encore de r\'{e}sultats.

\subsection{Principaux r\'{e}sultats}

Si on d\'{e}finit $A_{k}=\underset{m}{\oplus }(E_{k,m}T_{X}^{\ast })_{x}$
l'alg\`{e}bre des op\'{e}rateurs diff\'{e}rentiels en un point $x\in X,$
celle-ci peut-\^{e}tre vue comme une repr\'{e}sentation du groupe
lin\'{e}aire $Gl_{n}.$ On sait alors que l'on a une d\'{e}composition de
cette repr\'{e}sentation en somme directe de repr\'{e}sentations
irr\'{e}ductibles de Schur. Ainsi Demailly \cite{De95} a caract\'{e}ris\'{e}
les fibr\'{e}s de jets d'ordre 2, de degr\'{e} $m$ :
\begin{equation*}
Gr^{\bullet }E_{2,m}T_{X}^{\ast }=\underset{\lambda _{1}+2\lambda _{2}=m}{%
\oplus }\Gamma ^{(\lambda _{1},\lambda _{2},0)}T_{X}^{\ast },
\end{equation*}
o\`{u} $\Gamma $ est le foncteur de Schur.

\textit{\bigskip }

Le premier r\'{e}sultat de cet article est donn\'{e} par l'\'{e}tude
alg\'{e}brique de $A_{3},$ par la th\'{e}orie classique des invariants. On
obtient une caract\'{e}risation des jets d'ordre 3, en dimension 3:

\bigskip

\noindent \textbf{Th\'{e}or\`{e}me 1. }\textit{En dimension 3:}
\begin{equation*}
A_{3}=\mathbb{C[}f_{i}^{\prime },w_{ij},w_{ij}^{k},W],\text{ }1\leq i<j\leq
3,1\leq k\leq 3
\end{equation*}
\noindent \textit{o\`{u}} $W=$ $\left|
\begin{array}{ccc}
f_{1}^{\prime } & f_{2}^{\prime } & f_{3}^{\prime } \\
f_{1}^{\prime \prime } & f_{2}^{\prime \prime } & f_{3}^{\prime \prime } \\
f_{1}^{\prime \prime \prime } & f_{2}^{\prime \prime \prime } &
f_{3}^{\prime \prime \prime }
\end{array}
\right| ,w_{ij}=f_{i}^{\prime }f_{j}^{\prime \prime }-f_{i}^{\prime \prime
}f_{j}^{\prime },$

\noindent $w_{ij}^{k}=(f_{k}^{\prime })^{4}d(\frac{w_{ij}}{(f_{k}^{\prime
})^{3}})=f_{k}^{\prime }(f_{i}^{\prime }f_{j}^{\prime \prime \prime
}-f_{i}^{\prime \prime \prime }f_{j}^{\prime })-3f_{k}^{\prime \prime
}(f_{i}^{\prime }f_{j}^{\prime \prime }-f_{i}^{\prime \prime }f_{j}^{\prime
}).$

\textit{De plus, }$\deg .tr(\mathbb{C}(f_{i}^{\prime
},w_{ij},w_{ij}^{k},W))=7$\textit{\ (le calcul de l'id\'{e}al des relations
entre les g\'{e}n\'{e}rateurs est fait dans \cite{Rou04}).}

\bigskip

\noindent Cette \'{e}tude alg\'{e}brique a conduit \`{a} des applications
g\'{e}om\'{e}triques au niveau des fibr\'{e}s de jets.

\noindent Le r\'{e}sultat principal est la caract\'{e}risation du gradu\'{e}
du fibr\'{e} des jets d'ordre 3 en dimension 3:

\bigskip

\noindent \textbf{Th\'{e}or\`{e}me 2. }\textit{Soit X une vari\'{e}t\'{e}
complexe de dimension 3, alors:}
\begin{equation*}
Gr^{\bullet }E_{3,m}T_{X}^{\ast }=\underset{0\leq \gamma \leq \frac{m}{5}}{%
\oplus }(\underset{\{\lambda _{1}+2\lambda _{2}+3\lambda _{3}=m-\gamma ;%
\text{ }\lambda _{i}-\lambda _{j}\geq \gamma ,\text{ }i<j\}}{\oplus }\Gamma
^{(\lambda _{1},\lambda _{2},\lambda _{3})}T_{X}^{\ast })
\end{equation*}
\noindent \textit{o\`{u}} $\Gamma $ \textit{est le foncteur de Schur.}

\bigskip

\noindent Un calcul de type Riemann-Roch fournit alors:

\bigskip

\noindent \textbf{Proposition }\textit{Soit X une hypersurface lisse de
degr\'{e} }$d$ de $\mathbb{P}^{4}$, alors
\begin{equation*}
\chi (X,E_{3,m}T_{X}^{\ast })=\frac{m^{9}}{81648\times 10^{6}}%
d(389d^{3}-20739d^{2}+185559d-358873)+O(m^{8})
\end{equation*}

\noindent \textbf{Corollaire }\textit{Pour} $d\geq 43,$ $\chi
(X,E_{3,m}T_{X}^{\ast })\sim \alpha (d)m^{9}$ \textit{avec} $\alpha (d)>0.$

\section{Pr\'{e}liminaires}

Cette section a pour but de rappeler la construction des espaces de jets de
Demailly, les bases de la th\'{e}orie de la repr\'{e}sentation du groupe
lin\'{e}aire $Gl_{n}\mathbb{C}$ et celles de la th\'{e}orie classique des
invariants qui seront utilis\'{e}es de mani\`{e}re cruciale dans la preuve
des th\'{e}or\`{e}mes 1 et 2.

\subsection{Espaces des jets}

Soit $X$ une vari\'{e}t\'{e} complexe de dimension $n$. On d\'{e}finit le
fibr\'{e} $J_{k}\rightarrow X$ des $k$-jets de germes de courbes dans $X$,
comme \'{e}tant l'ensemble des classes d'\'{e}quivalence des applications
holomorphes $f:(\mathbb{C},0)\rightarrow (X,x)$ modulo la relation
d'\'{e}quivalence suivante: $f\sim g$ si et seulement si toutes les
d\'{e}riv\'{e}es $f^{(j)}(0)=g^{(j)}(0)$ co\"{i}ncident pour $0\leq j\leq k$%
. L'application projection $J_{k}\rightarrow X$ est simplement $f\rightarrow
f(0)$. Gr\^{a}ce \`{a} la formule de Taylor appliqu\'{e}e \`{a} un germe $f$
au voisinage d'un point $x\in X,$ on peut identifier $J_{k,x}$ \`{a}
l'ensemble des $k-$uplets de vecteurs $(f^{\prime }(0),...,f^{(k)}(0))\in
\mathbb{C}^{nk}.$ Ainsi, $J_{k}$ est un fibr\'{e} holomorphe sur $X$ de
fibre $\mathbb{C}^{nk}.$ On peut voir qu'il ne s'agit pas d'un fibr\'{e}
vectoriel pour $k\geq 2$ (pour $k=1,$ c'est simplement le fibr\'{e} tangent $%
T_{X}$)$.$

\begin{definition}
\textit{Soit (X,V) une vari\'{e}t\'{e} dirig\'{e}e. Le fibr\'{e} }$%
J_{k}V\rightarrow X$\textit{\ est l'espace des }$k-$\textit{jets de courbes }%
$f:(\mathbb{C},0)\rightarrow X$\textit{\ tangentes \`{a} V, c'est-\`{a}-dire
telles que }$f^{\prime }(t)\in V_{f(t)}$\textit{\ pour }$t$\textit{\ au
voisinage de }$0,$\textit{\ l'application projection sur X \'{e}tant }$%
f\rightarrow f(0).$
\end{definition}

Nous pr\'{e}sentons la construction des espaces de jets introduits par J.-P.
Demailly dans \cite{De95}.

Soit $(X,V)$ une vari\'{e}t\'{e} dirig\'{e}e. On d\'{e}finit $(X^{\prime
},V^{\prime })$ par :

\noindent i) $X^{\prime }=P(V)$

\noindent ii) $V^{\prime }\subset T_{X^{\prime }}$ est le sous-fibr\'{e} tel
que pour chaque point $(x,[v])\in X^{\prime }$ associ\'{e} \`{a} un vecteur $%
v\in V_{x}\backslash \{0\}$ on a :

\noindent $V_{(x,[v])}^{\prime }=\{\xi \in T_{X^{\prime }};\pi _{\ast }\xi
\in \mathbb{C}v\}$ o\`{u} $\pi :X^{\prime }\rightarrow X$ est la projection
naturelle et $\pi _{\ast }:T_{X^{\prime }}\rightarrow \pi ^{\ast }T_{X}$

\noindent On a donc $V^{\prime }=\pi _{\ast }^{-1}(O_{X^{\prime }}(-1)).$

On d\'{e}finit par r\'{e}currence le fibr\'{e} de k-jets projectivis\'{e} $%
P_{k}V=X_{k}$ et le sous-fibr\'{e} associ\'{e} $V_{k}\subset T_{X_{k}}$ par:

\noindent $(X_{0},V_{0})=(X,V),(X_{k},V_{k})=(X_{k-1}^{\prime
},V_{k-1}^{\prime }).$ On a par construction:

\begin{equation*}
\dim X_{k}=n+k(r-1),rangV_{k}=r:=rangV
\end{equation*}

Soit $\pi _{k}$ la projection naturelle $\pi _{k}:X_{k}\rightarrow X_{k-1},$
on notera $\pi _{j,k}:X_{k}\rightarrow X_{j}$ la composition $\pi
_{j+1}\circ \pi _{j+2}\circ ...\circ \pi _{k},$ pour $j\leq k.$

Par d\'{e}finition, il y a une injection canonique $O_{P_{k}V}(-1)%
\hookrightarrow \pi _{k}^{\ast }V_{k-1}$ et on obtient un morphisme de
fibr\'{e}s en droites
\begin{equation*}
O_{P_{k}V}(-1)\rightarrow \pi _{k}^{\ast }V_{k-1}\overset{(\pi _{k})^{\ast
}(\pi _{k-1})_{\ast }}{\rightarrow }\ \pi _{k}^{\ast }O_{P_{k-1}V}(-1)
\end{equation*}
\noindent qui admet
\begin{equation*}
D_{k}=P(T_{P_{k-1}V/P_{k-2}V})\subset P_{k}V
\end{equation*}
\noindent comme diviseur de z\'{e}ros

\noindent Ainsi, on a:
\begin{equation*}
O_{P_{k}V}(1)=\pi _{k}^{\ast }O_{P_{k-1}V}(1)\otimes O(D_{k}).
\end{equation*}

\begin{remark}
Chaque application non constante $f:\Delta _{R}\rightarrow X$ de $(X,V)$ se
rel\`{e}ve en $f_{[k]}:\Delta _{R}\rightarrow P_{k}V$. En effet:

\noindent si $f$ n'est pas constante, on peut d\'{e}finir la tangente $%
[f^{\prime }(t)]$ (aux points stationnaires $f^{\prime
}(t)=(t-t_{0})^{s}u(t),[f^{\prime }(t_{0})]=[u(t_{0})]$ ) et $%
f_{[1]}(t)=(f(t),[f^{\prime }(t)]).$
\end{remark}

\subsection{Op\'{e}rateurs diff\'{e}rentiels sur les jets}

D'apr\`{e}s \cite{GG80}, on introduit le fibr\'{e} vectoriel des jets de
diff\'{e}rentielles, d'ordre $k$ et de degr\'{e} $m$, $E_{k,m}^{GG}V^{\ast
}\rightarrow X$ dont les fibres sont les polyn\^{o}mes \`{a} valeurs
complexes $Q(f^{\prime },f^{\prime \prime },...,f^{(k)})$ sur les fibres de $%
J_{k}V,$ de poids $m$ par rapport \`{a} l'action de $\mathbb{C}^{\ast }$:
\begin{equation*}
Q(\lambda f^{\prime },\lambda ^{2}f^{\prime \prime },...,\lambda
^{k}f^{(k)})=\lambda ^{m}Q(f^{\prime },f^{\prime \prime },...,f^{(k)})
\end{equation*}
pour tout $\lambda \in \mathbb{C}^{\ast }$ et $(f^{\prime },f^{\prime \prime
},...,f^{(k)})\in J_{k}V.$

$E_{k,m}^{GG}V^{\ast }$ admet une filtration canonique dont les termes
gradu\'{e}s sont
\begin{equation*}
Gr^{l}(E_{k,m}^{GG}V^{\ast })=S^{l_{1}}V^{\ast }\otimes S^{l_{2}}V^{\ast
}\otimes ...\otimes S^{l_{k}}V^{\ast },
\end{equation*}
o\`{u} $l:=(l_{1},l_{2},...,l_{k})\in \mathbb{N}^{k}$ v\'{e}rifie $%
l_{1}+2l_{2}+...+kl_{k}=m.$ En effet, en consid\'{e}rant l'expression de
plus haut degr\'{e} en les $(f_{i}^{(k)})$ qui intervient dans l'expression
d'un polyn\^{o}me homog\`{e}ne de poids $m$, on obtient une filtration
intrins\`{e}que:
\begin{equation*}
E_{k-1,m}^{GG}V^{\ast }=S_{0}\subset S_{1}\subset ...\subset S_{\left[ \frac{%
m}{k}\right] }=E_{k,m}^{GG}V^{\ast }
\end{equation*}
o\`{u}
\begin{equation*}
S_{i}/S_{i-1}\simeq S^{i}V^{\ast }\otimes E_{k,m-ki}^{GG}V^{\ast }.
\end{equation*}
Par r\'{e}currence, on obtient bien une filtration dont les termes
gradu\'{e}s sont ceux annonc\'{e}s plus haut.

D'apr\`{e}s \cite{De95}, on d\'{e}finit le sous-fibr\'{e} $E_{k,m}V^{\ast
}\subset E_{k,m}^{GG}V^{\ast },$ appel\'{e} le fibr\'{e} des jets de
diff\'{e}rentielles invariants d'ordre $k$ et de degr\'{e} $m$, i.e :
\begin{equation*}
Q((f\circ \phi )^{\prime },(f\circ \phi )^{\prime \prime },...,(f\circ \phi
)^{(k)})=\phi ^{\prime }(0)^{m}Q(f^{\prime },f^{\prime \prime },...,f^{(k)})
\end{equation*}
pour tout $\phi \in G_{k}$ le groupe des germes de k-jets de
biholomorphismes de $(\mathbb{C},0).$ Pour $G_{k}^{\prime }$ le sous-groupe
de $G_{k}$ des germes $\phi $ tangents \`{a} l'identit\'{e} ($\phi ^{\prime
}(0)=1)$ on a $E_{k,m}V^{\ast }=(E_{k,m}^{GG}V^{\ast })^{G_{k}^{\prime }}.$

La filtration canonique sur $E_{k,m}^{GG}V^{\ast }$ induit une filtration
naturelle sur $E_{k,m}V^{\ast }$ dont les termes gradu\'{e}s sont
\begin{equation*}
\left( \underset{l_{1}+2l_{2}+...+kl_{k}=m}{\oplus }S^{l_{1}}V^{\ast
}\otimes S^{l_{2}}V^{\ast }\otimes ...\otimes S^{l_{k}}V^{\ast }\right)
^{G_{k}^{\prime }}.
\end{equation*}
Le lien entre ces espaces d'op\'{e}rateurs diff\'{e}rentiels et les espaces
de jets construits pr\'{e}c\'{e}demment est donn\'{e} par:

\begin{theorem}
\label{t12}\cite{De95} \textit{Supposons que V a un rang }$r\geq 2$\textit{.}

\textit{\noindent Soit }$\pi _{0,k}:P_{k}V\rightarrow X$\textit{, et }$%
J_{k}V^{reg}$\textit{\ le fibr\'{e} des k-jets r\'{e}guliers i.e }$f^{\prime
}(0)\neq 0.$

\textit{\noindent i) Le quotient }$J_{k}V^{reg}/G_{k}$\textit{\ a la
structure d'un fibr\'{e} localement trivial au-dessus de X, et il y a un
plongement holomorphe }$J_{k}V^{reg}/G_{k}\rightarrow P_{k}V,$\textit{\ qui
identifie }$J_{k}V^{reg}/G_{k}$\textit{\ avec }$P_{k}V^{reg}.$

\textit{\noindent ii) Le faisceau image direct }$(\pi _{0,k})_{\ast
}O_{P_{k}V}(m)$\textit{\ peut \^{e}tre identifi\'{e} avec le faisceau des
sections holomorphes de }$E_{k,m}V^{\ast }.$

\textit{\noindent iii) Pour tout }$m>0$\textit{, le lieu de base du
syst\`{e}me lin\'{e}aire }$\left| O_{P_{k}V}(m)\right| $\textit{\ est
\'{e}gal \`{a} }$P_{k}V^{sing}$\textit{. De plus, }$O_{P_{k}V}(1)$\textit{\
est relativement big (i.e pseudo-ample ) au-dessus de X.}
\end{theorem}

\subsection{Th\'{e}orie classique des invariants}

Soit $f$ un polyn\^{o}me dont les variables sont des vecteurs, i.e un
polyn\^{o}me en les coordonn\'{e}es des vecteurs, d'un espace vectoriel
fix\'{e} $V$. Pour tous vecteurs $s,t$ on note par $D_{s}^{t}f$ le
r\'{e}sultat de la diff\'{e}rentiation de $f$ par rapport \`{a} $s$ dans la
direction de $t,$ i.e:
\begin{equation*}
D_{s}^{t}f=\underset{i}{\sum }t_{i}\frac{\partial f}{\partial s_{i}}
\end{equation*}
o\`{u} les $s_{i},t_{i}$ sont les coordonn\'{e}es des vecteurs $s$ et $t$
respectivement.

Les op\'{e}rateurs de la forme $D_{s}^{t}$ sont appel\'{e}s op\'{e}rateurs
de polarisation. Ils commutent avec l'action du groupe $Gl(V)$ sur
l'alg\`{e}bre des polyn\^{o}mes.

Consid\'{e}rons la somme directe de m copies de $V$ munie de l'action
naturelle de $Gl(V)$ et de l'action de $Gl_{m}$ qui commute avec celle-ci,
i.e pour $A\in Gl_{m},$ $A.(x_{1},...,x_{m})=(x_{1},...,x_{m})A^{-1}.$ Ainsi
chaque vecteur $x_{j}$ est remplac\'{e} par une combinaison lin\'{e}aire de
vecteurs $x_{1},...,x_{m}$ avec des coefficients pris dans la j-\`{e}me
colonne de $A^{-1}.$ Cette action induit une action de $Gl_{m}$ sur
l'alg\`{e}bre des polyn\^{o}mes en les variables $x_{1},...,x_{m}.$
Explicitement, la matrice $A=(a_{ij})$ agit sur un polyn\^{o}me $f$ comme
suit:
\begin{equation*}
(Af)(x_{1},...,x_{m})=f(\underset{i}{\sum }a_{i1}x_{i},...,\underset{i}{\sum
}a_{im}x_{i}).
\end{equation*}
Si un polyn\^{o}me $f$ dont les variables sont des vecteurs a pour degr\'{e}
$p$ en la variable $x$, alors l'op\'{e}rateur $P_{x}=\frac{1}{p!}%
D_{x}^{x_{1}}...D_{x}^{x_{p}}$ (o\`{u} $x_{1},...,x_{p}$ n'apparaissent pas
dans l'expression de $f$) transforme $f$ en un polyn\^{o}me qui est
sym\'{e}trique et multi-lin\'{e}aire en $x_{1},...,x_{p}.$ On peut retrouver
$f$ \`{a} partir de $P_{x}f$ \ en substituant $x$ \`{a} la place de $%
x_{1},...,x_{p}.$ Si $f$ est homog\`{e}ne en toutes ses variables, si l'on
r\'{e}p\`{e}te l'op\'{e}ration pr\'{e}c\'{e}dente avec toutes les variables,
on obtient une forme multi-lin\'{e}aire $Pf$ appel\'{e}e \textit{la
polarisation compl\`{e}te} de $f.$ On retrouve $f$ en y substituant les
variables originelles.

\begin{definition}
(cf. \cite{Popov}) \textit{Soit }$F$\textit{\ une forme multi-lin\'{e}aire
en les variables }$u_{1},...,u_{l}$\textit{\ o\`{u} les }$u_{i}$\textit{\
sont des vecteurs d'un espace vectoriel }$V.$\textit{\ Soient }$%
x_{1},...,x_{m}$\textit{\ m vecteurs de }$V.$\textit{\ On d\'{e}finit }$%
S^{x_{1},...,x_{m}}(F),$\textit{\ l'espace vectoriel engendr\'{e} par tous
les polyn\^{o}mes obtenus en substituant les variables }$x_{1},...,x_{m}$%
\textit{\ aux variables }$u_{1},...,u_{l}$\textit{\ en permettant les
r\'{e}p\'{e}titions. Cet espace est clairement invariant sous l'action de }$%
Gl_{m}.$
\end{definition}

Soit $G$ un groupe lin\'{e}aire arbitraire agissant sur un espace vectoriel
de dimension $n.$ On consid\`{e}re le probl\`{e}me de trouver les $G-$%
invariants d'un syst\`{e}me de vecteurs de $V$, i.e les polyn\^{o}mes
invariants sous l'action de $G$ dans la somme directe de plusieurs copies de
$V.$ Il est clair que l'alg\`{e}bre de tous les $G-$invariants d'un
syst\`{e}me de vecteurs est lin\'{e}airement engendr\'{e} par les invariants
qui sont homog\`{e}nes en chaque variable. Si $f$ est un tel invariant, sa
polarisation compl\`{e}te en est un aussi. Ainsi si l'on est capable de
trouver tous les invariants multi-lin\'{e}aires, alors on obtient tous les
invariants homog\`{e}nes en y substituant de nouvelles variables (en
permettant les r\'{e}p\'{e}titions).

\begin{definition}
(cf.\cite{Popov}) \textit{Un ensemble }$\{F_{\alpha }\}$\textit{\ de formes
multi-lin\'{e}aires G-invariantes est appel\'{e} syst\`{e}me complet de
G-invariants d'un syst\`{e}me de m vecteurs si les espaces de polyn\^{o}mes }%
$S^{x_{1},...,x_{m}}(F)$\textit{\ associ\'{e}s aux formes }$F_{\alpha }$%
\textit{\ engendrent l'alg\`{e}bre de tous les G-invariants du syst\`{e}me
de vecteurs }$x_{1},...,x_{m}.$
\end{definition}

\begin{theorem}
(\cite{Popov})\label{t5} \textit{Soit V un espace vectoriel de dimension n.}

\noindent\textit{1) Tout syst\`{e}me complet de G-invariants d'un
syst\`{e}me de n vecteurs est aussi un syst\`{e}me complet pour tout nombre
de vecteurs.}

\noindent\textit{2) Si }$G\subset SL(V)$\textit{\ alors tout syst\`{e}me
complet de G-invariants d'un syst\`{e}me de }$n-1$\textit{\ vecteurs auquel
on ajoute la forme ''det'' est un syst\`{e}me complet de }$G$\textit{%
-invariants pour tout nombre de vecteurs.}
\end{theorem}

On rappelle qu'un groupe $G$ est dit lin\'{e}airement r\'{e}ductif si tout $%
G $-module V de dimension finie est semi-simple. On a alors le
th\'{e}or\`{e}me de Hilbert:

\begin{theorem}
(cf.\cite{Popov}) \textit{Soit }$G\subset Gl(V)$\textit{\ un groupe
r\'{e}ductif. Alors il existe un syst\`{e}me fini complet de }$G$\textit{%
-invariants.}
\end{theorem}

Dans le cas des groupes qui ne sont pas r\'{e}ductifs il y a quelques
r\'{e}sultats connus et des conjectures \`{a} propos du 14$^{\text{\`{e}me}}$
probl\`{e}me de Hilbert sur l'existence d'un syst\`{e}me fini de
g\'{e}n\'{e}rateurs de l'alg\`{e}bre des invariants. Le cas g\'{e}n\'{e}ral
se ram\`{e}ne au cas des groupes unipotents. Nagata (1959) a construit un
exemple de groupe unipotent dont l'alg\`{e}bre des invariants n'a pas de
syst\`{e}me fini de g\'{e}n\'{e}rateurs. Les r\'{e}sultats positifs
d\'{e}coulent du

\begin{theorem}
(cf.\cite{Popov}) (Principe de Grosshans) \textit{Soit }$G$\textit{\ un
groupe alg\'{e}brique qui agit rationnellement sur une }$k$\textit{%
-alg\`{e}bre A, et }$H$\textit{\ un sous-groupe ferm\'{e} de }$G$\textit{.
Alors:}
\begin{equation*}
A^{H}\cong (k[G]^{H}\otimes A)^{G}.
\end{equation*}
\end{theorem}

\noindent Si G est r\'{e}ductif et A de type fini, cela ram\`{e}ne le
probl\`{e}me de savoir si $A^{H}$ est de type fini \`{a} celui de savoir si $%
k[G]^{H}=k[G/H]$ est de type fini. D'o\`{u} la d\'{e}finition suivante:

\begin{definition}
(cf.\cite{Popov}) \textit{Un sous-groupe }$H$\textit{\ d'un groupe
r\'{e}ductif }$G$\textit{\ est appel\'{e} sous-groupe de Grosshans s'il
v\'{e}rifie les conditions: }$H$\textit{\ est ferm\'{e}, }$G/H$\textit{\ est
quasi-affine, }$k[G/H]$\textit{\ est de type fini.}
\end{definition}

On peut alors substituer au probl\`{e}me de Hilbert le probl\`{e}me suivant
propos\'{e} par K. Pommerening \cite{Pom87}: \textit{Trouver les
sous-groupes de Grosshans de }$Gl_{n}$\textit{\ ou plus g\'{e}n\'{e}ralement
d'un groupe r\'{e}ductif }$G.$

\noindent On a alors la conjecture de Popov-Pommerening \cite{Pom87}:

\begin{conjecture}
\label{co1}\textit{Tout sous-groupe unipotent r\'{e}gulier, i.e
normalis\'{e} par un tore maximal, d'un groupe r\'{e}ductif est de Grosshans.%
}
\end{conjecture}

L. Tan \cite{Tan89} a montr\'{e} que cette conjecture est vraie pour tous
les sous-groupes de $Gl_{n}(k),Sl_{n}(k),PSl_{n}(k)$ ($k$ corps
alg\'{e}briquement clos) pour $n\leq 5.$

\subsection{Th\'{e}orie de la repr\'{e}sentation}

Cette partie rappelle bri\`{e}vement la th\'{e}orie de la repr\'{e}sentation
de $Gl(V),$ o\`{u} $V$ est un espace vectoriel complexe de dimension finie $%
r $.

\bigskip

A l'ensemble des r-uplets d\'{e}croissants $(a_{1},...,a_{r})\in \mathbb{Z}%
^{r},a_{1}\geq a_{2}...\geq a_{r},$ on associe de mani\`{e}re fonctorielle
une collection d'espaces vectoriels $\Gamma ^{(a_{1},...,a_{r})}V$ qui
fournit la liste de toutes les repr\'{e}sentations polyn\^{o}miales
irr\'{e}ductibles du groupe lin\'{e}aire $Gl(V),$ \`{a} isomorphisme
pr\`{e}s. $\Gamma ^{\bullet }$ est appel\'{e} foncteur de Schur. Donnons une
description simple de ces foncteurs. Soit $\mathbb{U}_{r}=\left\{ \left(
\begin{array}{cc}
1 & \ast \\
0 & 1
\end{array}
\right) \right\} $ le groupe des matrices unipotentes triangulaires
sup\'{e}rieures $r\times r.$ Si tous les $a_{j}$ sont positifs, on
d\'{e}finit
\begin{equation*}
\Gamma ^{(a_{1},...,a_{r})}V\subset S^{a_{1}}V\otimes ...\otimes S^{a_{r}}V
\end{equation*}
comme \'{e}tant l'ensemble des polyn\^{o}mes $P(x_{1},...,x_{r})$ sur $%
(V^{\ast })^{r\text{ }}$qui sont homog\`{e}nes de degr\'{e} $a_{j}$ par
rapport \`{a} $x_{j}$ et qui sont invariants sous l'action \`{a} droite de $%
\mathbb{U}_{r}$ sur $(V^{\ast })^{r\text{ }}$ i.e tels que
\begin{equation*}
P(x_{1},...,x_{j-1},x_{j}+x_{k},x_{j+1},...,x_{r})=P(x_{1},...,x_{r})\text{
\ }\forall k<j.
\end{equation*}
Si $(a_{1},...,a_{r})$ n'est pas d\'{e}croissant alors on pose $\Gamma
^{(a_{1},...,a_{r})}V=0.$ Comme cas particuliers on retrouve les puissances
sym\'{e}triques et les puissances ext\'{e}rieures:
\begin{eqnarray*}
S^{k}V &=&\Gamma ^{(k,0,...,0)}V, \\
\wedge ^{k}V &=&\Gamma ^{(1,...,1,0,...,0)}V\text{ (avec }k\text{ indices 1),%
} \\
\det V &=&\Gamma ^{(1,...,1)}V.
\end{eqnarray*}
Les foncteurs de Schur satisfont la formule
\begin{equation*}
\Gamma ^{(a_{1}+l,...,a_{r}+l)}V=\Gamma ^{(a_{1},...,a_{r})}V\otimes (\det
V)^{l}
\end{equation*}
qui peut \^{e}tre utilis\'{e}e pour d\'{e}finir $\Gamma
^{(a_{1},...,a_{r})}V $ si l'on a des $a_{i}$ n\'{e}gatifs.

\bigskip

On fixe une base de V et on identifie $G=Gl(V)$ avec $Gl_{r}(\mathbb{C)}$.
On note $T=\{(x=diag(x_{1},...,x_{r})\}\subset G$ le sous-groupe des
matrices diagonales.

\begin{definition}
(cf.\cite{Fu.}) Un vecteur $e$ d'une repr\'{e}sentation E est appel\'{e}
vecteur de poids $\alpha =(\alpha _{1},...,\alpha _{r})$ (o\`{u} les $\alpha
_{i}$ sont des entiers) si
\begin{equation*}
x.e=x_{1}^{\alpha _{1}}...x_{r}^{\alpha _{r}}e\text{ pour tout }x\text{ de }%
T.
\end{equation*}
\end{definition}

\begin{proposition}
(cf.\cite{Fu.}) Toute repr\'{e}sentation E est somme directe de ses espaces
de poids:
\begin{equation*}
E=\oplus E_{\alpha },\text{ }E_{\alpha }=\{e\in E:x.e=x_{1}^{\alpha
_{1}}...x_{m}^{\alpha _{m}}e\text{ ,}\forall \text{ }x\in \text{ }T\text{ }%
\}.
\end{equation*}
\end{proposition}

\begin{definition}
\label{p9}(cf.\cite{Fu.}) Soit $B\subset G$ le groupe de Borel des matrices
triangulaires sup\'{e}rieures. Un vecteur $e$ d'une repr\'{e}sentation E est
appel\'{e} \textit{vecteur de plus haut poids} si $B.e=\mathbb{C}^{\ast }.e.$
\end{definition}

\begin{proposition}
\label{p6}(cf.\cite{Fu.}) Une repr\'{e}sentation (de dimension finie,
polyn\^{o}miale) E de $Gl_{r}(\mathbb{C)}$ est irr\'{e}ductible si et
seulement si elle a un unique vecteur de plus haut poids, \`{a}
multiplication par un scalaire pr\`{e}s. De plus, deux repr\'{e}sentations
sont isomorphes si et seulement si leurs vecteurs de plus haut poids ont le
m\^{e}me poids.
\end{proposition}

Nous utiliserons aussi la semi-simplicit\'{e} des repr\'{e}sentations
holomorphes de $Gl_{r}(\mathbb{C)}:$

\begin{proposition}
\label{p7}(cf.\cite{Fu.}) Toute repr\'{e}sentation holomorphe de $Gl_{r}(%
\mathbb{C)}$ est somme directe de repr\'{e}sentations irr\'{e}ductibles.
\end{proposition}

Ainsi pour d\'{e}terminer compl\`{e}tement une repr\'{e}sentation holomorphe
de $Gl_{r}(\mathbb{C)}$, il suffit de d\'{e}terminer ses vecteurs de plus
haut poids.

\section{Etude alg\'{e}brique}

On d\'{e}finit: $A_{k}=\underset{m}{\oplus }(E_{k,m}T_{X}^{\ast })_{x}$
l'alg\`{e}bre des op\'{e}rateurs diff\'{e}rentiels en un point $x\in X.$

Soit $G_{k}^{^{\prime }}$ le groupe des reparam\'{e}trisations $\phi
(t)=t+b_{2}t^{2}+...+b_{k}t^{k}+O(t^{k+1})$ tangentes \`{a} l'identit\'{e}. $%
G_{k}^{^{\prime }}$ agit sur $(f^{\prime },f^{\prime \prime },...,f^{(k)})$
par action unipotente. Par exemple pour $k=3$, on a l'action:
\begin{equation*}
(f\circ \phi )^{\prime }=f^{\prime };(f\circ \phi )^{\prime \prime
}=f^{\prime \prime }+2b_{2}f^{\prime };(f\circ \phi )^{\prime \prime \prime
}=f^{\prime \prime \prime }+6b_{2}f^{\prime \prime }+6b_{3}f^{\prime }
\end{equation*}

\noindent Donc une repr\'{e}sentation:
\begin{equation*}
G_{3}^{^{\prime }}\hookrightarrow U(3):\phi \rightarrow \left(
\begin{array}{ccc}
1 & 0 & 0 \\
2b_{2} & 1 & 0 \\
6b_{3} & 6b_{2} & 1
\end{array}
\right)
\end{equation*}

\noindent D\'{e}terminer $A_{k}$ revient donc \`{a} d\'{e}terminer $(\mathbb{%
C}[(f^{\prime }),(f^{\prime \prime }),...,(f^{(k)})])^{G_{k}^{^{\prime }}}.$

\bigskip

\noindent En dimension 2, on a $G_{2}^{^{\prime }}=U(2).$ Les invariants par
le groupe unipotent sont bien connus (cf.\cite{Procesi}). Ainsi:
\begin{eqnarray*}
(i)\text{ }A_{1} &=&\mathbb{C}[f_{1}^{\prime },f_{2}^{\prime }], \\
(ii)\text{ }A_{2} &=&\mathbb{C}[f_{1}^{\prime },f_{2}^{\prime },w_{12}]\text{
o\`{u} }w_{12}=f_{1}^{\prime }f_{2}^{\prime \prime }-f_{1}^{\prime \prime
}f_{2}^{\prime }.
\end{eqnarray*}

\bigskip

\noindent On a la propri\'{e}t\'{e} suivante:

\begin{proposition}
\label{p5}
\begin{equation*}
A_{n}=A_{n}[f_{1}^{\prime -1}]\cap A_{n}[f_{2}^{\prime -1}]
\end{equation*}
\end{proposition}

\begin{proof}
Il suffit de prouver $A_{n}[f_{1}^{\prime -1}]\cap A_{n}[f_{2}^{\prime
-1}]\subset $ $A_{n}.$ Soit $F\in A_{n}[f_{1}^{\prime -1}]\cap
A_{n}[f_{2}^{\prime -1}]:F=\frac{P}{(f_{1}^{\prime })^{l}}=\frac{Q}{%
(f_{2}^{\prime })^{m}}.$ Ainsi: $(f_{2}^{\prime })^{m}P=(f_{1}^{\prime
})^{l}Q$ et $(f_{1}^{\prime })^{l}$ divise P, donc $F\in \mathbb{C[}%
f^{\prime },f^{\prime \prime },f^{\prime \prime \prime }].$ De plus F est
invariant par reparam\'{e}trisation \ donc $F\in $ $A_{n}.$
\end{proof}

\bigskip

\subsection{Etude de la dimension 3 et preuve du th\'{e}or\`{e}me 1}

Nous \'{e}tudions maintenant la dimension 3:

\noindent $G_{3}^{^{\prime }}=\left\{ \left(
\begin{array}{ccc}
1 & 0 & 0 \\
2b_{2} & 1 & 0 \\
6b_{3} & 6b_{2} & 1
\end{array}
\right) \right\} \subset U(3).$

\noindent Faisons le lien avec la th\'{e}orie classique des invariants
(partie 3 des pr\'{e}liminaires).

\noindent $G_{3}^{^{\prime }}$ agit sur $\left(
\begin{array}{ccc}
f_{1}^{\prime } & f_{2}^{\prime } & f_{3}^{\prime } \\
f_{1}^{\prime \prime } & f_{2}^{\prime \prime } & f_{3}^{\prime \prime } \\
f_{1}^{\prime \prime \prime } & f_{2}^{\prime \prime \prime } &
f_{3}^{\prime \prime \prime }
\end{array}
\right) $ par multiplication \`{a} gauche.

\noindent Consid\'{e}rons l'action de $GL_{3}:$%
\begin{equation*}
A\in GL_{3},A.\left(
\begin{array}{ccc}
f_{1}^{\prime } & f_{2}^{\prime } & f_{3}^{\prime } \\
f_{1}^{\prime \prime } & f_{2}^{\prime \prime } & f_{3}^{\prime \prime } \\
f_{1}^{\prime \prime \prime } & f_{2}^{\prime \prime \prime } &
f_{3}^{\prime \prime \prime }
\end{array}
\right) =\left(
\begin{array}{ccc}
f_{1}^{\prime } & f_{2}^{\prime } & f_{3}^{\prime } \\
f_{1}^{\prime \prime } & f_{2}^{\prime \prime } & f_{3}^{\prime \prime } \\
f_{1}^{\prime \prime \prime } & f_{2}^{\prime \prime \prime } &
f_{3}^{\prime \prime \prime }
\end{array}
\right) A^{-1}.
\end{equation*}

\noindent Cette action induit une action sur les polyn\^{o}mes $P(f^{\prime
},f^{\prime \prime },f^{\prime \prime \prime })$ qui commute avec celle de $%
G_{3}^{^{\prime }}.$ Ainsi on a une action de $GL_{3}$ qui laisse $A_{3}$
invariant.

\noindent Nous cherchons \`{a} d\'{e}terminer les invariants par $%
G_{3}^{^{\prime }}$ du syst\`{e}me de vecteurs $(x_{1},x_{2},x_{3})$ o\`{u} $%
x_{i}=\left(
\begin{array}{c}
f_{i}^{\prime } \\
f_{i}^{\prime \prime } \\
f_{i}^{\prime \prime \prime }
\end{array}
\right) .$

Appliquons le th\'{e}or\`{e}me \ref{t5} des pr\'{e}liminaires \`{a} notre
situation. On a bien $G_{3}^{^{\prime }}\subset SL_{3}.$ Il nous suffit donc
de connaitre un syst\`{e}me complet de $G_{3}^{^{\prime }}$-invariants pour
deux vecteurs i.e en dimension 2. Cela nous est donn\'{e} par le
th\'{e}or\`{e}me annonc\'{e} par J.P. Demailly dont nous donnons ici une
d\'{e}monstration:

\begin{theorem}
\label{t8}(Demailly) \textit{En dimension 2:}
\begin{equation*}
A_{3}=\mathbb{C[}f_{1}^{\prime },f_{2}^{\prime
},w_{12}^{1},w_{12}^{2}][w_{12}]
\end{equation*}

\noindent \textit{o\`{u} }$w_{12}^{i}=(f_{i}^{\prime })^{4}d(\frac{w_{12}}{%
(f_{i}^{\prime })^{3}})=f_{i}^{\prime }(f_{1}^{\prime }f_{2}^{\prime \prime
\prime }-f_{1}^{\prime \prime \prime }f_{2}^{\prime })-3f_{i}^{^{\prime
\prime }}(f_{1}^{\prime }f_{2}^{\prime \prime }-f_{1}^{\prime \prime
}f_{2}^{\prime })$

\noindent \textit{et} $(\mathcal{R)}:$ $3(w_{12})^{2}=f_{2}^{\prime
}w_{12}^{1}-f_{1}^{\prime }w_{12}^{2}.$
\end{theorem}

\noindent La d\'{e}monstration n\'{e}cessite deux lemmes:

\begin{lemma}
\label{le1}$w_{12}$ \textit{est quadratique sur} $\mathbb{C[}f_{1}^{\prime
},f_{2}^{\prime },w_{12}^{2},w_{12}^{1}].$
\end{lemma}

\begin{proof}
\noindent Par $(\mathcal{R)},$ $w_{12}$ est alg\'{e}brique sur $\mathbb{C[}%
f_{1}^{\prime },f_{2}^{\prime },w_{12}^{2},w_{12}^{1}]$ de degr\'{e} 2 ou 1.

\noindent Supposons qu'il existe deux polyn\^{o}mes P et Q tels
que:\noindent
\begin{equation*}
P(f_{1}^{\prime },f_{2}^{\prime
},w_{12}^{2},w_{12}^{1})w_{12}=Q(f_{1}^{\prime },f_{2}^{\prime
},w_{12}^{2},w_{12}^{1}).
\end{equation*}

\noindent Par $(\mathcal{R)}$ on remplace $w_{12}^{2}$ par $\frac{%
f_{2}^{\prime }w_{12}^{1}-3(w_{12})^{2}}{f_{1}^{\prime }}$ dans P et Q.

\noindent Ainsi on obtient une \'{e}galit\'{e}, apr\`{e}s multiplication par
$(f_{1}^{\prime })^{m}$ avec $m$ suffisamment grand, entre deux
polyn\^{o}mes en les variables $\{f_{1}^{\prime },f_{2}^{\prime
},w_{12},w_{12}^{1}\}$ qui sont alg\'{e}briquement libres. Mais l'un des
polyn\^{o}mes a toutes ses puissances en $w_{12}$ impaires et l'autre,
paires; ce qui implique $P=Q=0.$

\noindent Ainsi le degr\'{e} de $w_{12}$ est 2.
\end{proof}

\begin{lemma}
\label{le2}$\{f_{1}^{\prime },f_{2}^{\prime },w_{12}^{2},w_{12}^{1}\}$
\textit{sont alg\'{e}briquement libres.}
\end{lemma}

\begin{proof}
$w_{12}$ est alg\'{e}brique sur $\mathbb{C(}f_{1}^{\prime },f_{2}^{\prime
},w_{12}^{2},w_{12}^{1})$ donc
\begin{eqnarray*}
\deg .tr(\mathbb{C(}f_{1}^{\prime },f_{2}^{\prime },w_{12}^{2},w_{12}^{1}))
&=&\deg .tr.(\mathbb{C(}f_{1}^{\prime },f_{2}^{\prime
},w_{12}^{2},w_{12}^{1},w_{12})) \\
&\geq &\deg .tr(\mathbb{C(}f_{1}^{\prime },f_{2}^{\prime
},w_{12},w_{12}^{1}))=4.
\end{eqnarray*}
\end{proof}

\bigskip

\noindent On peut maintenant passer \`{a} la d\'{e}monstration du
th\'{e}or\`{e}me \ref{t8}:

\bigskip

\begin{proof}
D'apr\`{e}s la proposition \ref{p5} on est ramen\'{e} \`{a} d\'{e}terminer $%
A_{3}[f_{1}^{\prime -1}]\cap A_{3}[f_{2}^{\prime -1}].$ On consid\`{e}re la
reparam\'{e}trisation $\phi =f_{1}^{-1}$ sur la carte $(f_{1}^{\prime }\neq
0).$ Soit $P\in A_{3}.$ Donc $P(f\circ \phi )=(\phi ^{\prime })^{m}P(f)\circ
\phi .$ Remarquons maintenant par le calcul:
\begin{eqnarray*}
(f_{2}\circ f_{1}^{-1})^{\prime } &=&\frac{f_{2}^{\prime }}{f_{1}^{^{\prime
}}}\circ f_{1}^{-1}, \\
(f_{2}\circ f_{1}^{-1})^{\prime \prime } &=&\frac{w_{12}}{(f_{1}^{\prime
})^{3}}\circ f_{1}^{-1}, \\
(f_{2}\circ f_{1}^{-1})^{\prime \prime \prime } &=&\frac{w_{12}^{1}}{%
(f_{1}^{\prime })^{5}}\circ f_{1}^{-1}.
\end{eqnarray*}

\noindent Ainsi $P\in \mathbb{C[}f_{1}^{\prime },f_{2}^{\prime
},w_{12},w_{12}^{1}][f_{1}^{\prime -1}]$ et donc $A_{3}[f_{1}^{\prime -1}]=%
\mathbb{C[}f_{1}^{\prime },f_{2}^{\prime },w_{12},w_{12}^{1}][f_{1}^{\prime
-1}].$ Par sym\'{e}trie: $A_{3}[f_{2}^{\prime -1}]=\mathbb{C[}f_{1}^{\prime
},f_{2}^{\prime },w_{12},w_{12}^{2}][f_{2}^{\prime -1}].$

\noindent L'inclusion
\begin{equation*}
\mathbb{C[}f_{1}^{\prime },f_{2}^{\prime
},w_{12},w_{12}^{1},w_{12}^{2}]\subset \mathbb{C[}f_{1}^{\prime
},f_{2}^{\prime },w_{12},w_{12}^{1}][f_{1}^{\prime -1}]\cap \mathbb{C[}%
f_{1}^{\prime },f_{2}^{\prime },w_{12},w_{12}^{2}][f_{2}^{\prime -1}]
\end{equation*}
est imm\'{e}diate puisque par $(\mathcal{R)}$:
\begin{equation*}
w_{12}^{2}\in \mathbb{C[}f_{1}^{\prime },f_{2}^{\prime
},w_{12},w_{12}^{1}][f_{1}^{\prime -1}]\text{ et }w_{12}^{1}\in \mathbb{C[}%
f_{1}^{\prime },f_{2}^{\prime },w_{12},w_{12}^{2}][f_{2}^{\prime -1}].
\end{equation*}

\noindent Il reste donc \`{a} montrer
\begin{equation*}
\mathbb{C[}f_{1}^{\prime },f_{2}^{\prime },w_{12},w_{12}^{1}][f_{1}^{\prime
-1}]\cap \mathbb{C[}f_{1}^{\prime },f_{2}^{\prime
},w_{12},w_{12}^{2}][f_{2}^{\prime -1}]\subset \mathbb{C[}f_{1}^{\prime
},f_{2}^{\prime },w_{12},w_{12}^{1},w_{12}^{2}].
\end{equation*}

\noindent Soit $F\in \mathbb{C[}f_{1}^{\prime },f_{2}^{\prime
},w_{12},w_{12}^{1}][f_{1}^{\prime -1}]\cap \mathbb{C[}f_{1}^{\prime
},f_{2}^{\prime },w_{12},w_{12}^{2}][f_{2}^{\prime -1}]:$%
\begin{equation*}
F=\frac{P(f_{1}^{\prime };f_{2}^{\prime };w_{12};w_{12}^{1})}{(f_{1}^{\prime
})^{l}}=\frac{Q(f_{1}^{\prime };f_{2}^{\prime };w_{12};w_{12}^{2})}{%
(f_{2}^{\prime })^{m}}
\end{equation*}
\noindent Par $(\mathcal{R)}$:
\begin{eqnarray*}
P(f_{1}^{\prime };f_{2}^{\prime };w_{12};w_{12}^{1}) &=&P_{1}(f_{1}^{\prime
};f_{2}^{\prime };w_{12}^{1};w_{12}^{2})w_{12}+P_{2}(f_{1}^{\prime
};f_{2}^{\prime };w_{12}^{1};w_{12}^{2}), \\
Q(f_{1}^{\prime };f_{2}^{\prime };w_{12};w_{12}^{2}) &=&Q_{1}(f_{1}^{\prime
};f_{2}^{\prime };w_{12}^{1};w_{12}^{2})w_{12}+Q_{2}(f_{1}^{\prime
};f_{2}^{\prime };w_{12}^{1};w_{12}^{2})
\end{eqnarray*}

\noindent Ainsi:
\begin{eqnarray*}
((f_{2}^{\prime })^{m}P_{1}(f_{1}^{\prime };f_{2}^{\prime
};w_{12}^{1};w_{12}^{2})-(f_{1}^{\prime })^{l}Q_{1}(f_{1}^{\prime
};f_{2}^{\prime };w_{12}^{1};w_{12}^{2}))w_{12} \\
+((f_{2}^{\prime })^{m}P_{2}(f_{1}^{\prime };f_{2}^{\prime
};w_{12}^{1};w_{12}^{2})-(f_{1}^{\prime })^{l}Q_{2}(f_{1}^{\prime
};f_{2}^{\prime };w_{12}^{1};w_{12}^{2})) &=&0.
\end{eqnarray*}
\noindent Or $w_{12}$ est quadratique sur $\mathbb{C[}f_{1}^{\prime
},f_{2}^{\prime },w_{12}^{2},w_{12}^{1}]$ donc:
\begin{equation*}
(f_{2}^{\prime })^{m}P_{i}(f_{1}^{\prime };f_{2}^{\prime
};w_{12}^{1};w_{12}^{2})-(f_{1}^{\prime })^{l}Q_{i}(f_{1}^{\prime
};f_{2}^{\prime };w_{12}^{1};w_{12}^{2})=0,\ i=1,2.
\end{equation*}
\noindent $\{f_{1}^{\prime },f_{2}^{\prime },w_{12}^{2},w_{12}^{1}\}$ sont
alg\'{e}briquement libres donc:
\begin{equation*}
P_{i}(f_{1}^{\prime };f_{2}^{\prime };w_{12}^{1};w_{12}^{2})=(f_{1}^{\prime
})^{l}R_{i}(f_{1}^{\prime };f_{2}^{\prime };w_{12}^{1};w_{12}^{2}).
\end{equation*}
\noindent Et le r\'{e}sultat est prouv\'{e}.
\end{proof}

\bigskip

\noindent On peut maintenant caract\'{e}riser les op\'{e}rateurs
diff\'{e}rentiels d'ordre 3 en dimension 3.

\noindent En notant $u_{i}=\left(
\begin{array}{c}
u_{i}^{1} \\
u_{i}^{2} \\
u_{i}^{3}
\end{array}
\right) $ et en d\'{e}finissant:
\begin{eqnarray*}
F_{1}(u_{1}) &=&u_{1}^{1}; \\
F_{2}(u_{1},u_{2}) &=&u_{1}^{1}u_{2}^{2}-u_{1}^{2}u_{2}^{1}; \\
F_{3}(u_{1},u_{2},u_{3})
&=&u_{3}^{1}(u_{1}^{1}u_{2}^{3}-u_{1}^{3}u_{2}^{1})-3u_{3}^{2}(u_{1}^{1}u_{2}^{2}-u_{1}^{2}u_{2}^{1}).
\end{eqnarray*}
on obtient que l'ensemble $\{F_{1},F_{2},F_{3}\}$ de formes
multilin\'{e}aires $G_{3}^{^{\prime }}$-invariantes est un syst\`{e}me
complet de $G_{3}^{^{\prime }}$-invariants d'un syst\`{e}me de 2 vecteurs.

\noindent Par application du th\'{e}or\`{e}me \ref{t5} de Popov, on obtient
la preuve du th\'{e}or\`{e}me 1 et donc, la caract\'{e}risation
alg\'{e}brique de l'alg\`{e}bre $A_{3}$ des germes d'op\'{e}rateurs
invariants en dimension 3:

\bigskip

\begin{proof}
Il ne reste qu'\`{a} justifier l'assertion sur le degr\'{e} de
transcendance. Mais celle-ci est une cons\'{e}quence imm\'{e}diate du
th\'{e}or\`{e}me \ref{t12} qui identifie $E_{k,m}T_{X,x}^{\ast }$ avec les
sections de $O_{P_{k}V}(m)$ au-dessus de $(\pi _{0,k})^{-1}(x).$
\end{proof}

\begin{remark}
\noindent 1) Pour tout $k,$ $G_{k}^{^{\prime }}\subset SL_{k}$, donc par le
raisonnement pr\'{e}c\'{e}dent pour d\'{e}terminer $A_{k}$ en toute
dimension il suffit de d\'{e}terminer $A_{k}$ en dimension $k-1.$

\noindent 2) On a montr\'{e} que le groupe $G_{3}^{^{\prime }}=\left\{
\left(
\begin{array}{ccc}
1 & 0 & 0 \\
2b_{2} & 1 & 0 \\
6b_{3} & 6b_{2} & 1
\end{array}
\right) \right\} \subset U(3)$ est un groupe de Grosshans de $GL_{3}$ i.e $%
\mathbb{C[}GL_{3}]^{G_{3}^{^{\prime }}}$ est une alg\`{e}bre de type fini.
De plus, ce groupe n'est pas r\'{e}gulier i.e normalis\'{e} par un tore
maximal car:
\begin{eqnarray*}
\left(
\begin{array}{ccc}
\lambda _{1} & 0 & 0 \\
0 & \lambda _{2} & 0 \\
0 & 0 & \lambda _{3}
\end{array}
\right) \left(
\begin{array}{ccc}
1 & 0 & 0 \\
2b_{2} & 1 & 0 \\
6b_{3} & 6b_{2} & 1
\end{array}
\right) \left(
\begin{array}{ccc}
\lambda _{1}^{-1} & 0 & 0 \\
0 & \lambda _{2}^{-1} & 0 \\
0 & 0 & \lambda _{3}^{-1}
\end{array}
\right) &=& \\
\left(
\begin{array}{ccc}
1 & 0 & 0 \\
2b_{2}\lambda _{1}^{-1}\lambda _{2} & 1 & 0 \\
6b_{3}\lambda _{1}^{-1}\lambda _{3} & 6b_{2}\lambda _{2}^{-1}\lambda _{3} & 1
\end{array}
\right) &\notin &G_{3}^{^{\prime }}.
\end{eqnarray*}
On ne peut donc pas appliquer le r\'{e}sultat de L. Tan \cite{Tan89} sur la
conjecture \ref{co1} de Popov-Pommerening cit\'{e} dans les
pr\'{e}liminaires pour montrer que $G_{3}^{^{\prime }}$ est un sous-groupe
de Grosshans.

\noindent 3) Sans l'utilisation du th\'{e}or\`{e}me de Popov, la
d\'{e}termination par un calcul ''\`{a} la main'' des g\'{e}n\'{e}rateurs de
$A_{3} $ semble difficile.
\end{remark}

\section{Applications g\'{e}om\'{e}triques et preuve du th\'{e}o\-r\`{e}me 2}

Il s'agit d'\'{e}tudier le fibr\'{e} $E_{3,m}T_{X}^{\ast }$ en dimension 3
pour obtenir sa filtration en repr\'{e}sentations irr\'{e}ductibles de Schur
qui nous permettra, par un calcul de Riemann-Roch, de calculer sa
caract\'{e}ristique d'Euler. Rappelons (cf. introduction) que $%
E_{3,m}T_{X}^{\ast }$ est muni d'une filtration dont les termes gradu\'{e}s
sont
\begin{equation*}
Gr^{\bullet }E_{3,m}T_{X}^{\ast }=\left( \underset{l_{1}+2l_{2}+3l_{3}=m}{%
\oplus }S^{l_{1}}T_{X}^{\ast }\otimes S^{l_{2}}T_{X}^{\ast }\otimes
S^{l_{3}}T_{X}^{\ast }\right) ^{G_{3}^{\prime }}.
\end{equation*}
D'apr\`{e}s la th\'{e}orie de la repr\'{e}sentation, ces termes gradu\'{e}s
se d\'{e}composent en repr\'{e}sentations irr\'{e}ductibles de $%
Gl(T_{X}^{\ast })$: les repr\'{e}sentations de Schur. La caract\'{e}risation
alg\'{e}brique pr\'{e}c\'{e}dente va nous permettre de trouver les
repr\'{e}sentations irr\'{e}ductibles qui interviennent dans cette
d\'{e}composition.

\bigskip

\noindent \noindent Pour cela, on a besoin de la filtration des 3-jets en
dimension 2:

\begin{theorem}
\label{t9}\textit{En dimension 2 on a:}
\begin{equation*}
Gr^{\bullet }E_{3,m}T_{X}^{\ast }=\underset{0\leq \gamma \leq \frac{m}{5}}{%
\oplus }(\underset{\{\lambda _{1}+2\lambda _{2}=m-\gamma ;\text{ }\lambda
_{1}-\lambda _{2}\geq \gamma ;\text{ }\lambda _{2}\geq \gamma \}}{\oplus }%
\Gamma ^{(\lambda _{1},\lambda _{2})}T_{X}^{\ast })
\end{equation*}
\end{theorem}

\begin{proof}
On sait que
\begin{equation*}
A_{3}=\mathbb{C[}f_{1}^{\prime },f_{2}^{\prime
},w_{12}^{1},w_{12}^{2}][w_{12}]
\end{equation*}
\noindent o\`{u} $w_{12}^{i}=(f_{i}^{\prime })^{4}d(\frac{w_{12}}{%
(f_{i}^{\prime })^{3}})=f_{i}^{\prime }(f_{1}^{\prime }f_{2}^{\prime \prime
\prime }-f_{1}^{\prime \prime \prime }f_{2}^{\prime })-3f_{i}^{\prime \prime
}(f_{1}^{\prime }f_{2}^{\prime \prime }-f_{1}^{\prime \prime }f_{2}^{\prime
})$

\noindent et $3(w_{12})^{2}=f_{2}^{\prime }w_{12}^{1}-f_{1}^{\prime
}w_{12}^{2}.$

\noindent $A_{3,m}$ est une repr\'{e}sentation polyn\^{o}miale de $GL_{2}.$
La th\'{e}orie de la repr\'{e}sentation (proposition \ref{p6} et \ref{p7})
nous dit que $A_{3,m}$ est somme directe de repr\'{e}sentations
irr\'{e}ductibles qui sont d\'{e}termin\'{e}es par les vecteurs de plus haut
poids.

\noindent Rappelons (d\'{e}finition \ref{p9}) qu'un vecteur est vecteur de
plus haut poids s'il est invariant sous l'action de $U(2)=\left\{ \left(
\begin{array}{cc}
1 & \ast \\
0 & 1
\end{array}
\right) \right\} .$

\noindent Ici:
\begin{equation*}
V=\{(f_{1}^{\prime })^{\alpha }(w_{12}^{1})^{\gamma }(w_{12})^{\beta }\text{
}/\text{ }\alpha +5\gamma +3\beta =m\}
\end{equation*}
est clairement un ensemble de vecteurs de plus haut poids, de poids
\begin{equation*}
(\alpha +\beta +2\gamma ,\beta +\gamma ).
\end{equation*}

\noindent On en d\'{e}duit que chaque repr\'{e}sentation $\Gamma ^{(\lambda
_{1},\lambda _{2})}$ v\'{e}rifiant
\begin{equation*}
\{\lambda _{1}+2\lambda _{2}=m-\gamma ;\lambda _{1}-\lambda _{2}\geq \gamma
;\lambda _{2}\geq \gamma \}
\end{equation*}
appara\^{i}t une et une seule fois dans les repr\'{e}sentations
d\'{e}termin\'{e}es par cet ensemble de vecteurs de plus haut poids. En
effet, soit $(\lambda _{1},\lambda _{2})$ un tel couple alors
\begin{equation*}
\{\alpha =\lambda _{1}-\lambda _{2}-\gamma ;\beta =\lambda _{2}-\gamma \}
\end{equation*}
et $(\alpha ,\beta ,\gamma )$ sont d\'{e}termin\'{e}s de mani\`{e}re unique.

\noindent On a donc:
\begin{equation*}
Gr^{\bullet }E_{3,m}T_{X}^{\ast }\supset \underset{0\leq \gamma \leq \frac{m%
}{5}}{\oplus }(\underset{\{\lambda _{1}+2\lambda _{2}=m-\gamma ;\text{ }%
\lambda _{1}-\lambda _{2}\geq \gamma ;\text{ }\lambda _{2}\geq \gamma \}}{%
\oplus }\Gamma ^{(\lambda _{1},\lambda _{2})}T_{X}^{\ast }).
\end{equation*}

\noindent Pour avoir l'\'{e}galit\'{e} il suffit de montrer que l'ensemble V
est l'ensemble de tous les vecteurs de plus haut poids, i.e:
\begin{equation*}
V=(A_{3,m})^{U(2)}.
\end{equation*}
\noindent Soit $P\in (A_{3,m})^{U(2)}:P=P_{1}+P_{2}.w_{12},$ avec $P_{i}\in
\mathbb{C[}f_{1}^{\prime },f_{2}^{\prime },w_{12}^{1},w_{12}^{2}].$

\noindent Soit $u\in U(2):u.P=u.P_{1}+(u.P_{2}).w_{12}$ car $%
u.w_{12}=w_{12}. $

\noindent Donc $u.P=P\Leftrightarrow u.P_{i}=P_{i}$ (car $w_{12}$ est
quadratique par le lemme \ref{le1}).

\noindent Donc pour d\'{e}terminer $(A_{3,m})^{U(2)},$ il nous suffit de
d\'{e}terminer $\mathbb{C[}f_{1}^{\prime },f_{2}^{\prime
},w_{12}^{1},w_{12}^{2}]^{U(2)}.$

\noindent Soit: $u=\left(
\begin{array}{cc}
1 & \lambda \\
0 & 1
\end{array}
\right) \in U(2):$

\noindent On a les relations suivantes:
\begin{eqnarray*}
u.f_{1}^{\prime } &=&f_{1}^{\prime }; \\
u.f_{2}^{\prime } &=&\lambda f_{1}^{\prime }+f_{2}^{\prime }; \\
u.w_{12}^{1} &=&w_{12}^{1}; \\
u.w_{12}^{2} &=&w_{12}^{2}+\lambda w_{12}^{1}.
\end{eqnarray*}

\noindent Rappelons que $\{f_{1}^{\prime },f_{2}^{\prime
},w_{12}^{2},w_{12}^{1}\}$ sont alg\'{e}briquement libres par le lemme \ref
{le2}, donc d\'{e}terminer $\mathbb{C[}f_{1}^{\prime },f_{2}^{\prime
},w_{12}^{1},w_{12}^{2}]^{U(2)}$ revient \`{a} d\'{e}terminer les invariants
du groupe unipotent $U(2)$ qui sont bien connus en th\'{e}orie classique des
invariants (cf.\cite{Procesi} p.87). Donc on a l'\'{e}galit\'{e}:
\begin{equation*}
\mathbb{C[}f_{1}^{\prime },f_{2}^{\prime },w_{12}^{1},w_{12}^{2}]^{U(2)}=%
\mathbb{C[}f_{1}^{\prime },w_{12}^{1},f_{2}^{\prime
}w_{12}^{1}-f_{1}^{\prime }w_{12}^{2}]=\mathbb{C[}f_{1}^{\prime
},w_{12}^{1},(w_{12})^{2}].
\end{equation*}

\noindent Finalement on obtient l'inclusion:
\begin{equation*}
(A_{3,m})^{U(2)}\subset \mathbb{C[}f_{1}^{\prime },w_{12}^{1},w_{12}].
\end{equation*}

\noindent Par l'unicit\'{e} de $(\alpha ,\beta ,\gamma )$ vue
pr\'{e}c\'{e}demment on obtient bien:
\begin{equation*}
(A_{3,m})^{U(2)}=V.
\end{equation*}
\end{proof}

\bigskip

\noindent On passe maintenant \`{a} la preuve du th\'{e}or\`{e}me 2:

\bigskip

\begin{proof}
On suit le m\^{e}me sch\'{e}ma que dans la preuve pr\'{e}c\'{e}dente.

\noindent Soit
\begin{equation*}
V=\{(f_{1}^{\prime })^{\alpha }(w_{12}^{1})^{\gamma }(w_{12})^{\beta
}W^{\delta }\text{ }/\text{ }\alpha +5\gamma +3\beta +6\delta =m\}.
\end{equation*}

\noindent V est un ensemble de vecteurs de plus haut poids de poids
\begin{equation*}
(\alpha +\beta +2\gamma +\delta ;\beta +\gamma +\delta ;\delta ).
\end{equation*}

\noindent Soit $(\lambda _{1},\lambda _{2},\lambda _{3})$ v\'{e}rifiant:
\begin{equation*}
(\mathcal{P}):\{\lambda _{1}+2\lambda _{2}+3\lambda _{3}=m-\gamma ,0\leq
\gamma \leq \frac{m}{5};\lambda _{i}-\lambda _{j}\geq \gamma ,i<j\}.
\end{equation*}

Comme pr\'{e}c\'{e}demment on obtient que chaque repr\'{e}sentation $\Gamma
^{(\lambda _{1},\lambda _{2},\lambda _{3})}T_{X}^{\ast }$ o\`{u} $(\lambda
_{1},\lambda _{2},\lambda _{3})$ v\'{e}rifie $(\mathcal{P})$ appara\^{i}t
une et une seule fois dans les repr\'{e}sentations d\'{e}termin\'{e}es par
cet ensemble de vecteurs de plus haut poids. En effet, soit $(\lambda
_{1},\lambda _{2},\lambda _{3})$ v\'{e}rifiant $(\mathcal{P}).$

\noindent Alors:
\begin{equation*}
\{\alpha =\lambda _{1}-\lambda _{2}-\gamma ;\text{ }\beta =\lambda
_{2}-\lambda _{3}-\gamma ;\text{ }\delta =\lambda _{3}\}
\end{equation*}
et $(\alpha ,\beta ,\gamma ,\delta )$ sont d\'{e}termin\'{e}s de mani\`{e}re
unique.

\noindent Donc on a l'inclusion:
\begin{equation*}
Gr^{\bullet }E_{3,m}T_{X}^{\ast }\supset \underset{0\leq \gamma \leq \frac{m%
}{5}}{\oplus }(\underset{\{\lambda _{1}+2\lambda _{2}+3\lambda _{3}=m-\gamma
;\text{ }\lambda _{i}-\lambda _{j}\geq \gamma ,\text{ }i<j\}}{\oplus }\Gamma
^{(\lambda _{1},\lambda _{2},\lambda _{3})}T_{X}^{\ast }).
\end{equation*}

\noindent Pour avoir l'\'{e}galit\'{e} il suffit \`{a} nouveau de montrer
que V est l'ensemble de tous les vecteurs de plus haut poids de $A_{3,m}$
i.e: $V=(A_{3,m})^{U(2)}.$

L'id\'{e}e importante ici est d'utiliser un argument qui apparait dans la
preuve du th\'{e}or\`{e}me \ref{t5} de Popov \cite{Popov} et permet de voir
que le r\'{e}sultat obtenu pour la dimension 2 implique le r\'{e}sultat pour
la dimension 3.

\noindent Si $(x_{1},x_{2},x_{3})$ est un syst\`{e}me de vecteurs en
position g\'{e}n\'{e}rale tel que
\begin{equation*}
\det (x_{1},x_{2},x_{3})=0
\end{equation*}
alors par l'action de $U(3)$ on se ram\`{e}ne au syst\`{e}me $%
(x_{1},x_{2},0).$

\noindent Soit $P\in (A_{3,m})^{U(3)}$, un vecteur de plus haut poids.
Montrons que
\begin{equation*}
P\in \mathbb{C[}f_{1}^{\prime },w_{12}^{1},w_{12},W]
\end{equation*}
par r\'{e}currence sur $m.$ Pour $m=0,$ c'est trivial.

\noindent Supposons maintenant $(A_{3,p})^{U(3)}\subset \mathbb{C[}%
f_{1}^{\prime },w_{12}^{1},w_{12},W]$ pour $p<m.$ Montrons que le
r\'{e}sultat est vrai pour $m.$ Consid\'{e}rons $P_{1}$ la restriction de $P$
\`{a} l'hypersurface $(W=0).$ Par l'invariance de $P_{1}$ sous l'action de $%
U(3)$ et la remarque pr\'{e}c\'{e}dente montrant que par $U(3)$ on
transforme le syst\`{e}me $(x_{1},x_{2},x_{3}),$ en position
g\'{e}n\'{e}rale, en le syst\`{e}me $(x_{1},x_{2},0)$, on obtient que $P_{1}$
ne d\'{e}pend que des deux premiers vecteurs i.e $P_{1}$ est un vecteur de
plus haut poids de dinension 2, donc par le th\'{e}or\`{e}me \ref{t9} $%
P_{1}\in \mathbb{C[}f_{1}^{\prime },w_{12}^{1},w_{12}].$

\noindent $P-P_{1}$ est un polyn\^{o}me qui s'annule sur l'hypersurface $%
(W=0).$ Par le Nullstellensatz, on obtient que $(P-P_{1})\in \sqrt{(W)}$
donc par l'irr\'{e}ductibilit\'{e} de $W$ on a:
\begin{equation*}
P=P_{1}+W.P_{2}.
\end{equation*}

\noindent Il est clair que $P_{2}\in (A_{3,m-6})^{U(3)}$ donc par
hypoth\`{e}se de r\'{e}currence
\begin{equation*}
P_{2}\in \mathbb{C[}f_{1}^{\prime },w_{12}^{1},w_{12},W]
\end{equation*}
\ et de m\^{e}me pour $P$.

\noindent On en d\'{e}duit que $(A_{3,m})^{U(3)}\subset \mathbb{C[}%
f_{1}^{\prime },w_{12}^{1},w_{12},W].$

\noindent Donc $V=(A_{3,m})^{U(3)}$ par l'unicit\'{e} de $(\alpha ,\beta
,\gamma ,\delta ).$

\noindent Le th\'{e}or\`{e}me est d\'{e}montr\'{e}.
\end{proof}

\section{Calculs de caract\'{e}ristiques d'Euler}

Soit $X\subset \mathbb{P}^{4}$ une hypersurface lisse et irr\'{e}ductible de
degr\'{e} $d$. Gr\^{a}ce aux filtrations obtenues dans la section
pr\'{e}c\'{e}dente nous allons pouvoir calculer les diff\'{e}rentes
caract\'{e}ristiques d'Euler qui nous int\'{e}ressent. Les calculs ont
\'{e}t\'{e} faits sur le logiciel Maple et d\'{e}taill\'{e}s dans \cite
{Rou04}.

\subsection{Calcul des classes de Chern}

Soit $c_{i}=c_{i}(T_{X}).$

\begin{proposition}
\begin{eqnarray*}
c_{1}^{3} &=&(5-d)^{3}d, \\
c_{1}c_{2} &=&d(5-d)(d^{2}-5d+10), \\
c_{3} &=&d(-d^{3}+5d^{2}-10d+10).
\end{eqnarray*}
\end{proposition}

\begin{proof}
On a la suite exacte du fibr\'{e} normal:
\begin{equation*}
0\rightarrow T_{X}\rightarrow T_{\mathbb{P}^{4}\left| X\right. }\rightarrow
\mathcal{O}_{X}(d)\rightarrow 0.
\end{equation*}
Donc par d\'{e}finition des classes de Chern:
\begin{equation*}
(1+c_{1}+c_{2}+c_{3})(1+dh)=(1+h)^{5}
\end{equation*}
o\`{u}: $h=c_{1}(\mathcal{O}_{X}(1)).$ Donc, par identification on obtient
les identit\'{e}s:
\begin{eqnarray*}
c_{1} &=&(5-d)h, \\
c_{2}+dc_{1}h &=&10h^{2}, \\
c_{3}+dc_{2}h &=&10h^{3}, \\
h^{3} &=&d.
\end{eqnarray*}
d'o\`{u} les relations annonc\'{e}es:
\begin{eqnarray*}
c_{1}^{3} &=&(5-d)^{3}d, \\
c_{3} &=&10d-10d^{2}+d^{3}(5-d), \\
c_{1}c_{2} &=&(5-d)10d-d^{2}(5-d)^{2}.
\end{eqnarray*}
\end{proof}

\subsection{Les 1-jets}

L'absence de 1-jets d\'{e}finis globalement est bien connue:

\begin{proposition}
\textbf{(\cite{Sa2}})\label{p11}
\begin{equation*}
H^{0}(X,S^{m}T_{X}^{\ast })=0\text{\textit{\ pour }}m\geq 1.
\end{equation*}
\end{proposition}

\noindent Donc les 1-jets ne pourront pas \^{e}tre utilis\'{e}s.

\begin{remark}
\label{r1}Un calcul de type Riemann-Roch donne:
\begin{equation*}
\chi (X,S^{m}T_{X}^{\ast })=\frac{m^{5}}{120}%
(-c_{1}^{3}+2c_{1}c_{2}-c_{3})+O(m^{4})=\frac{m^{5}}{120}5d(3d-7)+O(m^{4}).
\end{equation*}
\end{remark}

\subsection{Les 2-jets}

On a la filtration \cite{De95}:
\begin{equation*}
Gr^{\bullet }E_{2,m}T_{X}^{\ast }=\underset{l_{1}+2l_{2}=m}{\oplus }\Gamma
^{(l_{1},l_{2},0)}T_{X}^{\ast }.
\end{equation*}
d'o\`{u} le calcul:

\begin{equation*}
\chi (X,E_{2,m}T_{X}^{\ast })=\chi (X,Gr^{\bullet }E_{2,m}T_{X}^{\ast })=%
\underset{l_{1}+2l_{2}=m}{\sum }\chi (X,\Gamma ^{(l_{1},l_{2},0)}T_{X}^{\ast
}).
\end{equation*}

\begin{proposition}
\begin{equation*}
\chi (X,E_{2,m}T_{X}^{\ast })=\frac{-m^{7}}{1837080}%
(89c_{1}^{3}-141c_{1}c_{2}+52c_{3})+O(m^{6})
\end{equation*}
\noindent Donc:
\begin{equation*}
\chi (X,E_{2,m}T_{X}^{\ast })=\frac{m^{7}}{1837080}%
(-5d(37d^{2}-452d+919))+O(m^{6})
\end{equation*}
\end{proposition}

On constate donc la n\'{e}gativit\'{e} de la caract\'{e}ristique d'Euler
pour $d$ suffisamment grand.

\begin{remark}
Pour les jets de Green-Griffiths :
\begin{equation*}
\chi (X,\mathcal{J}_{2,m})=\frac{m^{8}}{2^{3}8!}(\frac{-15}{8}%
c_{1}^{3}+3c_{1}c_{2}-\frac{9}{8}c_{3})+O(m^{7})
\end{equation*}
Donc:
\begin{equation*}
\chi (X,\mathcal{J}_{2,m})=\frac{-15m^{8}}{8^{3}7!}%
(d(51-25d+2d^{2})+O(m^{7}).
\end{equation*}
\end{remark}

\bigskip

\subsection{Les 3-jets}

Gr\^{a}ce \`{a} la filtration obtenue pr\'{e}c\'{e}demment dans le
th\'{e}or\`{e}me 2, on peut effectuer un calcul de Riemann-Roch :

\begin{proposition}
\begin{equation*}
\chi (X,E_{3,m}T_{X}^{\ast })=-m^{9}(\frac{43}{1417500000}c_{3}+\frac{29233}{%
408240000000}c_{1}^{3}-\frac{551}{5670000000}c_{1}c_{2})+O(m^{8}).
\end{equation*}
Donc:
\begin{equation*}
\chi (X,E_{3,m}T_{X}^{\ast })=\frac{m^{9}}{81648\times 10^{6}}%
d(389d^{3}-20739d^{2}+185559d-358873)+O(m^{8})
\end{equation*}
\end{proposition}

\begin{corollary}
\textit{Pour} $d\geq 43,$ $\chi (X,E_{3,m}T_{X}^{\ast })\sim \alpha (d)m^{9}$
avec $\alpha (d)>0.$
\end{corollary}

\begin{remark}
Pour les jets de Green-Griffiths on a :
\begin{eqnarray*}
\chi (X,\mathcal{J}_{3,m}) &=&\frac{-m^{11}}{6^{3}.11!}(\frac{575}{216}%
c_{1}^{3}-\frac{395}{108}c_{1}c_{2}+\frac{251}{216}c_{3})+O(m^{10}) \\
&=&\frac{m^{11}}{6^{3}.11!.216}d(36d^{3}-1980d^{2}+17985d-34885)+O(m^{10})
\end{eqnarray*}
Et on a la positivit\'{e} pour $d\geq 45.$
\end{remark}

\subsection{Le cas logarithmique}

Nous pouvons appliquer les r\'{e}sultats obtenus au cas logarithmique. Soit $%
X$ une vari\'{e}t\'{e} lisse complexe avec un diviseur \`{a} croisements
normaux $D$. En suivant \cite{Ii77}, on d\'{e}finit le faisceau cotangent
logarithmique
\begin{equation*}
\overline{T_{X}}^{\ast }=T_{X}^{\ast }(\log D)
\end{equation*}
comme le faisceau localement libre engendr\'{e} par $T_{X}^{\ast }$ et les
diff\'{e}rentielles logarithmiques $\frac{ds_{j}}{s_{j}},$ o\`{u} les $%
s_{j}=0$ sont les \'{e}quations locales des composantes irr\'{e}ductibles de
$D$.

Son dual, le fibr\'{e} tangent logarithmique
\begin{equation*}
\overline{T_{X}}=T_{X}(-\log D)
\end{equation*}
est le faisceau des germes de champs de vecteurs tangents \`{a} $D$.

De la m\^{e}me mani\`{e}re que dans le cas compact, on peut construire les
espaces de jets logarithmiques et les fibr\'{e}s d'op\'{e}rateurs
diff\'{e}rentiels associ\'{e}s $E_{k,m}\overline{T_{X}}^{\ast }$ (cf. \cite
{DL96}).

Soit $X\subset \mathbb{P}^{3}$ une surface lisse et irr\'{e}ductible de
degr\'{e} $d$. On consid\`{e}re la vari\'{e}t\'{e} logarithmique $(\mathbb{P}%
^{3},X).$

\subsubsection{Calcul des classes de Chern}

On pose :
\begin{equation*}
c_{i}=c_{i}(T_{\mathbb{P}^{3}}),\overline{c_{i}}=c_{i}(\overline{T_{\mathbb{P%
}^{3}}}).
\end{equation*}

\bigskip

Le calcul des classes de Chern est un peu plus long que dans le cas compact.

\begin{proposition}
\begin{eqnarray*}
\overline{c_{1}}^{3} &=&(4-d)^{3}, \\
\overline{c_{1}}\overline{c_{2}} &=&(4-d)(d^{2}-4d+6), \\
\overline{c_{3}} &=&-d^{3}+4d^{2}-6d+4.
\end{eqnarray*}
\end{proposition}

\begin{proof}
Pour la premi\`{e}re identit\'{e} il suffit de remarquer que :
\begin{equation*}
\overline{c_{1}}=-c_{1}(\overline{K_{\mathbb{P}^{3}}})=-c_{1}(\mathcal{O}_{%
\mathbb{P}^{3}}(d-4)).
\end{equation*}
La troisi\`{e}me vient du fait \cite{Ii78} que :
\begin{equation*}
\overline{c_{3}}=e(\mathbb{P}^{3}\backslash X)
\end{equation*}
o\`{u} $e$ d\'{e}signe la caract\'{e}ristique d'Euler.

\noindent On a $e(\mathbb{P}^{3}\backslash X)=e(\mathbb{P}^{3})-e(X)$ et $e(%
\mathbb{P}^{3})=4.$ Calculons $e(X)=c_{2}(T_{X}).$ Par la suite exacte
\begin{equation*}
0\rightarrow T_{X}\rightarrow T_{\mathbb{P}^{3}\left| X\right. }\rightarrow
\mathcal{O}_{X}(d)\rightarrow 0
\end{equation*}
on obtient
\begin{equation*}
c(T_{\mathbb{P}^{3}\left| X\right. })=c(T_{X})c(\mathcal{O}_{X}(d))
\end{equation*}
donc
\begin{equation*}
(1+h)^{4}=(1+c_{1}(T_{X})+c_{2}(T_{X})).(1+dh)
\end{equation*}
o\`{u} $h=c_{1}(\mathcal{O}_{X}(1))$ et $h^{2}=d.$ De plus
\begin{equation*}
c_{1}(T_{X})=-c_{1}(K_{X})=(4-d)h
\end{equation*}
donc on a l'\'{e}galit\'{e}
\begin{equation*}
1+4h+6h^{2}=(1+(4-d)h+c_{2}(T_{X})).(1+dh).
\end{equation*}
On a donc l'identit\'{e}
\begin{equation*}
e(X)=c_{2}(T_{X})=6h^{2}-(4-d)dh^{2}=d(6+d^{2}-4d)
\end{equation*}
d'o\`{u} finalement
\begin{equation*}
e(\mathbb{P}^{3}\backslash X)=4-d(6+d^{2}-4d).
\end{equation*}
Montrons la deuxi\`{e}me. Rappelons que par Riemann-Roch, si E est un
fibr\'{e} vectoriel de rang e sur $\mathbb{P}^{3},$ avec des classes de
Chern not\'{e}es $d_{i},$ alors (\cite{KO70} p508):
\begin{equation*}
\chi (\mathbb{P}^{3},E)=\frac{1}{6}(d_{1}^{3}-3d_{1}d_{2}+3d_{3})+\frac{1}{4}%
c_{1}(d_{1}^{2}-2d_{2})+\frac{1}{12}(c_{1}^{2}+c_{2})d_{1}+\frac{e}{24}%
c_{1}c_{2}.
\end{equation*}
Par Riemann-Roch:
\begin{equation*}
\chi (\overline{T_{\mathbb{P}^{3}}}^{\ast })=\frac{1}{6}(-\overline{c_{1}}%
^{3}+3\overline{c_{1}}\overline{c_{2}}-3\overline{c_{3}})+\frac{1}{4}c_{1}(%
\overline{c_{1}}^{2}-2\overline{c_{2}})-\frac{1}{12}(c_{1}^{2}+c_{2})%
\overline{c_{1}}+\frac{1}{8}c_{1}c_{2}.
\end{equation*}
Donc:
\begin{equation*}
\frac{1}{2}(\overline{c_{1}}\overline{c_{2}}-c_{1}\overline{c_{2}})=\chi (%
\overline{T_{\mathbb{P}^{3}}}^{\ast })-(\frac{1}{6}(-\overline{c_{1}}^{3}-3%
\overline{c_{3}})+\frac{1}{4}c_{1}\overline{c_{1}}^{2}-\frac{1}{12}%
(c_{1}^{2}+c_{2})\overline{c_{1}}+\frac{1}{8}c_{1}c_{2}).
\end{equation*}
Pour d\'{e}terminer $\overline{c_{1}}\overline{c_{2}},$ il suffit donc de
d\'{e}terminer $\chi (\overline{T_{\mathbb{P}^{3}}}^{\ast }).$

\noindent On a la suite exacte:
\begin{equation*}
0\rightarrow T_{\mathbb{P}^{3}}^{\ast }\rightarrow \overline{T_{\mathbb{P}%
^{3}}}^{\ast }\rightarrow \mathcal{O}_{X}\rightarrow 0
\end{equation*}
o\`{u} la fl\`{e}che $\phi :\overline{T_{\mathbb{P}^{3}}}^{\ast }\rightarrow
\mathcal{O}_{X}$ est donn\'{e}e par:

\noindent soit $x\in \mathbb{P}^{3},w\in (\overline{T_{\mathbb{P}^{3}}}%
^{\ast })_{x}$. Si $x\notin X:\phi (w)=0.$ Si $x\in X,w=f_{1}\frac{dz_{1}}{%
z_{1}}+f_{2}dz_{2}+f_{3}dz_{3},$ o\`{u} $(z_{1}=0)$ est une \'{e}quation
locale de X: $\phi (w)=f_{1}.$

\noindent Donc:
\begin{equation*}
\chi (\overline{T_{\mathbb{P}^{3}}}^{\ast })=\chi (\mathcal{O}_{X})+\chi (T_{%
\mathbb{P}^{3}}^{\ast }).
\end{equation*}
Calculons $\chi (T_{\mathbb{P}^{3}}^{\ast }):$

\noindent on a
\begin{equation*}
c_{1}(T_{\mathbb{P}^{3}}^{\ast })=c_{1}(\mathcal{O}_{\mathbb{P}%
^{3}}(-4))=-4\omega
\end{equation*}
o\`{u} $\omega $ est la classe d'un hyperplan, $c_{2}(T_{\mathbb{P}%
^{3}}^{\ast })=6\omega ^{2}$, $c_{3}(T_{\mathbb{P}^{3}}^{\ast })=-4.$

\noindent Par Riemann-Roch:
\begin{eqnarray*}
\chi (T_{\mathbb{P}^{3}}^{\ast }) &=&\frac{1}{6}((-4)^{3}+3\times 4\times
6-3\times 4)+\frac{1}{4}\times 4\times (4^{2}-2\times 6) \\
&&+\frac{1}{12}(4^{2}+6)\times (-4)+\frac{1}{8}\times 4\times 6 \\
&=&-1.
\end{eqnarray*}
Calculons $\chi (\mathcal{O}_{X}).$ On a la suite exacte:
\begin{equation*}
0\rightarrow \mathcal{O}(-X)\rightarrow \mathcal{O}_{\mathbb{P}%
^{3}}\rightarrow \mathcal{O}_{X}\rightarrow 0.
\end{equation*}
qui nous donne:
\begin{equation*}
\chi (\mathcal{O}_{X})=\chi (\mathcal{O}_{\mathbb{P}^{3}})-\chi (\mathcal{O}%
(-X)).
\end{equation*}
\noindent $\chi (\mathcal{O}_{\mathbb{P}^{3}})=\frac{1}{24}c_{1}c_{2}$ et $%
\chi (\mathcal{O}(-X))=\frac{-d^{3}}{6}+d^{2}-\frac{22}{12}d+\frac{1}{24}%
c_{1}c_{2}.$

\noindent Donc:
\begin{equation*}
\chi (\mathcal{O}_{X})=\frac{d^{3}}{6}-d^{2}+\frac{11}{6}d.
\end{equation*}
Finalement:
\begin{equation*}
\chi (\overline{T_{\mathbb{P}^{3}}}^{\ast })=\frac{d^{3}}{6}-d^{2}+\frac{11}{%
6}d-1.
\end{equation*}
$c_{1}=4\omega =\frac{4}{4-d}\overline{c_{1}}.$ Donc:
\begin{equation*}
\chi (\overline{T_{\mathbb{P}^{3}}}^{\ast })=\frac{1}{2}(1-\frac{4}{4-d})%
\overline{c_{1}}\overline{c_{2}}=\frac{-d}{2(4-d)}\overline{c_{1}}\overline{%
c_{2}}.
\end{equation*}
Et:
\begin{eqnarray*}
\overline{c_{1}}\overline{c_{2}} &=&\frac{2(d-4)}{d}(\chi (\overline{T_{%
\mathbb{P}^{3}}}^{\ast })-(\frac{1}{6}(-\overline{c_{1}}^{3}-3\overline{c_{3}%
})+\frac{1}{4}c_{1}\overline{c_{1}}^{2}-\frac{1}{12}(c_{1}^{2}+c_{2})%
\overline{c_{1}}+\frac{1}{8}c_{1}c_{2})) \\
&=&\frac{2(d-4)}{d}(\frac{d^{3}}{6}-d^{2}+\frac{11}{6}d-1-(\frac{1}{6}\times
((d-4)^{3}-3\times (-d^{3}+4d^{2}-6d+4)) \\
&&+\frac{1}{4}\times 4\times (4-d)^{2}-\frac{1}{12}\times 22\times (4-d)+3)
\\
&=&(4-d)(d^{2}-4d+6).
\end{eqnarray*}
\end{proof}

\subsubsection{Calcul des caract\'{e}risiques d'Euler}

Nous montrons d'abord que les filtrations restent les m\^{e}mes que dans le
cas compact mis-\`{a}-part que le fibr\'{e} tangent est remplac\'{e} par le
fibr\'{e} tangent logarithmique. Comme dans le cas compact (cf.\cite{De95}),
on munit le faisceau $\mathcal{O}(E_{k,m}^{GG}\overline{T_{X}}^{\ast })$\
des diff\'{e}rentielles de jets logarithmiques (i.e le faisceau localement
libre engendr\'{e} par tous les op\'{e}rateurs polyn\^{o}miaux en les d\'{e}%
riv\'{e}es d'ordre $1,2,...k$ de $f$, auxquelles on ajoute celles de la
fonction $\log (s_{j}(f))$\ le long de la j-\`{e}me composante de $D$) d'une
filtration dont les termes gradu\'{e}s sont
\begin{equation*}
Gr^{\bullet }E_{k,m}^{GG}\overline{T_{X}}^{\ast }=S^{l_{1}}\overline{T_{X}}%
^{\ast }\otimes ...\otimes S^{l_{k}}\overline{T_{X}}^{\ast },
\end{equation*}
$l:=(l_{1},l_{2},...,l_{k})\in \mathbb{N}^{k}$ v\'{e}rifie $%
l_{1}+2l_{2}+...+kl_{k}=m.$ Pour le faisceau des diff\'{e}rentielles de jets
invariants $\mathcal{O}(E_{k,m}\overline{T_{X}}^{\ast })$ (cf. \cite{DL96})
on obtient une filtration dont les termes gradu\'{e}s sont
\begin{equation*}
Gr^{\bullet }E_{k,m}\overline{T_{X}}^{\ast }=\left( \underset{%
l_{1}+2l_{2}+...+kl_{k}=m}{\oplus }S^{l_{1}}\overline{T_{X}}^{\ast }\otimes
...\otimes S^{l_{k}}\overline{T_{X}}^{\ast }\right) ^{G_{k}^{\prime }}
\end{equation*}
o\`{u} l'action de $G_{k}^{\prime }$ est \'{e}tendue de $U\backslash D,$ o%
\`{u} $U$ est un ouvert de $X,$ \`{a} $U$ gr\^{a}ce \`{a} l'isomorphisme
(cf. \cite{DL96}) $\overline{T_{X}}_{\left| U\right. }^{\ast }\rightarrow
T_{X,p}^{\ast }\times U$ , $(p\in U\backslash D).$ On a le diagramme
commutatif
\begin{equation*}
\begin{array}{ccc}
\scriptstyle\left( \underset{l_{1}+...+kl_{k}=m}{\oplus }S^{l_{1}}\overline{%
T_{X}}^{\ast }\otimes ...\otimes S^{l_{k}}\overline{T_{X}}^{\ast }\right)
_{\left| U\backslash D\right. }^{G_{k}^{\prime }} & {\small \rightarrow } & %
\scriptstyle\left( \underset{l_{1}+...+kl_{k}=m}{\oplus }S^{l_{1}}T_{X,p}^{%
\ast }\otimes ...\otimes S^{l_{k}}T_{X,p}^{\ast }\right) ^{G_{k}^{\prime }}%
{\small \times (U\backslash D)} \\
{\small \downarrow \scriptstyle Gl(}\scriptstyle\overline{T_{X}}^{\ast }%
{\small )} &  & {\small \downarrow \scriptstyle Gl}_{3}{\small \times id} \\
\scriptstyle\left( \underset{l_{1}+...+kl_{k}=m}{\oplus }S^{l_{1}}\overline{%
T_{X}}^{\ast }\otimes ...\otimes S^{l_{k}}\overline{T_{X}}^{\ast }\right)
_{\left| U\backslash D\right. }^{G_{k}^{\prime }} & {\small \rightarrow } & %
\scriptstyle\left( \underset{l_{1}+...+kl_{k}=m}{\oplus }S^{l_{1}}T_{X,p}^{%
\ast }\otimes ...\otimes S^{l_{k}}T_{X,p}^{\ast }\right) ^{G_{k}^{\prime }}%
{\small \times (U\backslash D)}
\end{array}
\end{equation*}

\noindent o\`{u} les fl\`{e}ches horizontales sont des isomorphismes qui
s'\'{e}tendent \`{a} $U.$ L'action de $Gl_{3}\times id$ s'\'{e}tend
clairement \`{a} $U,$ donc l'action de $Gl(\overline{T_{X}}^{\ast }),$ \`{a}
gauche dans le diagramme aussi. Les repr\'{e}sentations irr\'{e}ductibles
\`{a} droite s'identifient donc avec celles de gauche car les vecteurs de
plus haut poids et les poids s'identifient par la commutativit\'{e} du
diagramme. Ainsi on a bien la m\^{e}me d\'{e}composition en
repr\'{e}sentations irr\'{e}ductibles dans le cas logarithmique et dans le
cas compact i.e
\begin{eqnarray*}
Gr^{\bullet }E_{2,m}\overline{T_{X}}^{\ast } &=&\underset{l_{1}+2l_{2}=m}{%
\oplus }\Gamma ^{(l_{1},l_{2},0)}\overline{T_{X}}^{\ast }, \\
Gr^{\bullet }E_{3,m}\overline{T_{X}}^{\ast } &=&\underset{0\leq \gamma \leq
\frac{m}{5}}{\oplus }(\underset{\{\lambda _{1}+2\lambda _{2}+3\lambda
_{3}=m-\gamma ;\text{ }\lambda _{i}-\lambda _{j}\geq \gamma ,\text{ }i<j\}}{%
\oplus }\Gamma ^{(\lambda _{1},\lambda _{2},\lambda _{3})}\overline{T_{X}}%
^{\ast }).
\end{eqnarray*}

\bigskip

\noindent D'apr\`{e}s les calculs dans le cas compact, on obtient les
r\'{e}sultats suivants (calculs d\'{e}taill\'{e}s dans \cite{Rou04}):

\bigskip

\noindent Pour les 1-jets:
\begin{equation*}
\chi (S^{m}\overline{T_{\mathbb{P}^{3}}}^{\ast })=\frac{m^{5}}{120}%
(10d-20)+O(m^{4}).
\end{equation*}

\bigskip

\noindent Pour les 2-jets:
\begin{equation*}
\chi (E_{2,m}\overline{T_{\mathbb{P}^{3}}}^{\ast })=m^{7}(\frac{-37}{459270}%
d^{2}+\frac{247}{306180}d-\frac{1}{129})+O(m^{6}).
\end{equation*}

\bigskip

\noindent On obtient \`{a} nouveau la positivit\'{e} pour les 3-jets:
\begin{equation*}
\chi (E_{3,m}\overline{T_{\mathbb{P}^{3}}}^{\ast })=m^{9}(\frac{389}{%
81648000000}d^{3}-\frac{6913}{34020000000}d^{2}+\frac{6299}{4252500000}d-%
\frac{1513}{63787500})+O(m^{8}).
\end{equation*}

\begin{corollary}
\textit{Pour} $d\geq 34,$ $\chi (X,E_{3,m}\overline{T_{\mathbb{P}^{3}}}%
^{\ast })\sim \alpha (d)m^{9}$ avec $\alpha (d)>0.$
\end{corollary}

\bigskip

\bibliographystyle{amsplain}
\bibliography{jets}

\end{document}